\documentclass[11pt]{article}
\usepackage{amssymb}
\usepackage{amsmath}
\usepackage{times}

\title{Baire class one colorings\\ and a dichotomy for countable unions of $F_\sigma$ rectangles\indent}
\author{Dominique LECOMTE}
\date{March 2010}

\def\ufootnote#1{\let\savedthfn\thefootnote\let\thefootnote\relax
\footnote{#1}\let\thefootnote\savedthfn\addtocounter{footnote}{-1}}

\newcommand{\Ana}{{\it\Sigma}^{1}_{1}}

\newcommand{\Ca}{{\it\Pi}^{1}_{1}}
\newcommand{\Boraone}{{\it\Sigma}^{0}_{1}}

\newcommand{\Borel}{{\it\Delta}^{1}_{1}}

\newcommand{\borel}{{\bf\Delta}^{1}_{1}}

\newcommand{\boraone}{{\bf\Sigma}^{0}_{1}}
\newcommand{\boratwo}{{\bf\Sigma}^{0}_{2}}

\newcommand{\boraxi}{{\bf\Sigma}^{0}_{\xi}}

\newcommand{\borone}{{\bf\Delta}^{0}_{1}}
\newcommand{\bortwo}{{\bf\Delta}^{0}_{2}}

\newcommand{\bormone}{{\bf\Pi}^{0}_{1}}

\newcommand{\bormlxi}{{\bf\Pi}^{0}_{<\xi}}

\newcommand{\bormxi}{{\bf\Pi}^{0}_{\xi}}

\newcommand{\borxi}{{\bf\Delta}^{0}_{\xi}}
\newcommand{\borme}{{\bf\Pi}^{0}_{\eta}}

\newtheorem{thm} {Theorem} [section]
\newtheorem{defi} [thm] {Definition}

\newtheorem{lem} [thm] {Lemma}
\newtheorem{prop} [thm] {Proposition}

\begin{document}

\maketitle

\centerline{$\bullet$ Universit\' e Paris 6, Institut de Math\'ematiques de Jussieu, Projet Analyse Fonctionnelle}

\centerline{Tour 46-0, bo\^\i te 186, 4, place Jussieu, 75 252 Paris Cedex 05, France.}

\centerline{dominique.lecomte@upmc.fr}\bigskip

\centerline{$\bullet$ Universit\'e de Picardie, I.U.T. de l'Oise, site de Creil,}

\centerline{13, all\'ee de la fa\"\i encerie, 60 107 Creil, France.}
\bigskip\bigskip\bigskip\bigskip\bigskip\bigskip

\ufootnote{{\it 2010 Mathematics Subject Classification.}~Primary: 03E15, Secondary: 54H05}

\ufootnote{{\it Keywords and phrases.}~Baire class one, Borel chromatic number, Borel class, coloring,  dichotomy, Hurewicz, partition, product}

\noindent {\bf Abstract.} We study the Baire class one countable colorings, i.e., the countable partitions into $F_\sigma$ sets. Such a partition gives a covering of the diagonal into countably many $F_\sigma$ squares. This leads to the study of countable unions of $F_\sigma$ rectangles. We give a Hurewicz-like dichotomy for such countable unions.

\vfill\eject

\section{$\!\!\!\!\!\!$ Introduction}\indent

 The reader should see [K] for the standard descriptive set theoretic notation used in this paper. We study a definable coloring problem. We will need some more notation:\bigskip
 
\noindent\bf Notation.\rm ~The letters $X$, $Y$ will refer to some sets. We set 
$\Delta (X)\! :=\!\{ (x_0,x_1)\!\in\! X^2\mid x_0\! =\! x_1\}$.

\begin{defi} (1) Let $A\!\subseteq\! X^2$. We say that $A$ is a $digraph$ if 
$A\cap\Delta (X)\! =\!\emptyset$.\smallskip

\noindent (2) Let $A$ be a digraph. A $countable~coloring$ of $(X,A)$ is a map 
$c\! :\! X\!\rightarrow\!\omega$ such that $A$ does not meet $(c\!\times\! c)^{-1}\big(\Delta (\omega )\big)$.\end{defi}

 In [K-S-T], the authors characterize the analytic digraphs of having a Borel countable coloring. The characterization is given in terms of the following notion of comparison between relations.\bigskip

\noindent\bf Notation.\rm ~Let $X,Y$ be Polish spaces, $A$ (resp., $B$) a relation on $X$ (resp., $Y$), and $\bf\Gamma$ a class of sets.
$$(X,A)\preceq_{\bf\Gamma}(Y,B)~\Leftrightarrow ~\exists f\! :\! X\!\rightarrow\! Y~~{\bf\Gamma}\mbox{-measurable with }A\!\subseteq\! (f\!\times\! f)^{-1}(B).$$
In this case, we say that $f$ is a ${\bf\Gamma}$-$measurable\ homomorphism$ from $(X,A)$ into $(Y,B)$. This notion essentially makes sense for digraphs (we can take $f$ to be constant if $B$ is not a digraph).\bigskip

 We also have to introduce a minimum digraph without Borel countable coloring:\bigskip

\noindent $\bullet$ Let $\psi\! :\!\omega\!\rightarrow\! 2^{<\omega}$ be the natural bijection. More specifically, $\psi (0)\! :=\!\emptyset$ is the sequence of length $0$, $\psi (1)\! :=\! 0$, $\psi (2)\! :=\! 1$ are the sequences of length $1$, and so on. Note that $|\psi (n)|\!\leq\! n$ if $n\!\in\!\omega$. Let 
$n\!\in\!\omega$. As $|\psi (n)|\!\leq\! n$, we can define $s_n\! :=\!\psi (n)0^{n-|\psi (n)|}$. The crucial properties of the sequence $(s_{n})_{n\in\omega}$ are the following:\bigskip

- For each $s\!\in\! 2^{<\omega}$, there is $n\!\in\!\omega$ such that $s\!\subseteq\! s_n$ (we say that 
$(s_n)_{n\in\omega}$ is $dense$ in $2^{<\omega}$).\smallskip

- $|s_n|\! =\! n$.\bigskip

\noindent $\bullet$ We put
$\mathbb{G}_0\! :=\!\{ (s_n0\gamma ,s_n1\gamma )\mid n\!\in\!{\omega}\mbox{ and }
\gamma\!\in\! 2^{\omega}\}\!\subseteq\! 2^\omega\!\times\! 2^\omega$. Note that $\mathbb{G}_0$ is analytic since the map $(n,\gamma )\!\mapsto\! (s_n0\gamma ,s_n1\gamma )$ is continuous.\bigskip

 The previous definitions were given, when ${\bf\Gamma}\! =\!\borel$, in [K-S-T], where the following is proved:

\begin{thm} (Kechris, Solecki, Todor\v cevi\'c) Let $X$ be a Polish space, and $A$ an analytic relation on $X$. Then exactly one of the following holds:\smallskip  

\noindent (a) There is a Borel countable coloring of $(X,A)$, i.e., 
$(X,A)\preceq_{\borel}\big(\omega ,\neg\Delta (\omega )\big)$,\smallskip  

\noindent (b) $(2^{\omega},\mathbb{G}_0)\preceq_{\boraone}(X,A)$.\end{thm} 

 This result had several developments during the last years:\bigskip
 
\noindent - We can characterize the potentially closed sets via a Hurewicz-like test, and in finite dimension it is a consequence of the previous result. Let us specify this. The following definition can be found in [Lo2] (see Definition 3.3).
 
\begin{defi} (Louveau) Let $X,Y$ be Polish spaces, $A$ a Borel subset of $X\!\times\! Y$, and $\bf\Gamma$ a Borel class. We say that $A$ is $potentially~in~{\bf\Gamma}$ $\big($denoted 
$A\!\in\!\mbox{pot}({\bf\Gamma})\big)$ iff we can find a finer Polish topology $\sigma$ $($resp., $\tau )$ on $X$ $($resp., $Y)$ such that $A\!\in\! {\bf\Gamma}\big( (X, \sigma )\!\times\! (Y, \tau )\big)$.\end{defi}

 In particular, the potentially open sets are exactly the countable unions of Borel rectangles. A consequence of this is that the Borel hierarchy build on the Borel rectangles is exactly the hierarchy of the classes of the sets potentially in some Borel class.\bigskip
 
 The good notion of comparison to study the $\mbox{pot}({\bf\Gamma})$ sets is as follows. Let $X_0,X_1,Y_0,Y_1$ be Polish spaces, and  
$A^\varepsilon_0,A^\varepsilon_1$ disjoint analytic subsets of $X_\varepsilon\!\times\! Y_\varepsilon$. Then we set\bigskip

\leftline{$(X_0,Y_0,A^0_0,A^0_1)\leq (X_1,Y_1,A^1_0,A^1_1)\Leftrightarrow$}\smallskip

\rightline{$\exists f\! :\! X_0\!\rightarrow\! X_1~~\exists g\! :\! Y_0\!\rightarrow\! Y_1\mbox{ continuous with }
\forall\varepsilon\!\in\! 2~~A^0_\varepsilon\!\subseteq\! (f\!\times\! g)^{-1}(A^1_\varepsilon )\mbox{,}$}\bigskip

 The following theorem is proved in [L1], and is a consequence of Theorem 1.2:
 
\begin{thm} Let $X,Y$ be Polish spaces, and $A_0,A_1$ disjoint analytic subsets of $X\!\times\! Y$. Then exactly one of the following holds:\smallskip

\noindent (a) The set $A_0$ can be separated from $A_1$ by a $\mbox{pot}(\boraone )$ set,\smallskip

\noindent (b) $\big( 2^\omega ,2^\omega ,\Delta (2^\omega ),\mathbb{G}_0\big)\leq (X,Y,A_0,A_1)$.
\end{thm}
 
 In [L1], it is also proved that we cannot have $f$ one-to-one in Theorem 1.2.(b) in general. It is easy to check that Theorem 1.2 is also an easy consequence of Theorem 1.4. This means that the study of the Borel countable colorings is highly related to the study of countable unions of Borel rectangles.\bigskip
 
\noindent - We can extend Theorem 1.2 to any finite dimension, and also in infinite dimension if we change the space in which lives the infinite dimensional version of $\mathbb{G}_0$ (see [L2]).\bigskip
 
\noindent - B. Miller recently developped some techniques to recover many dichotomy results of descriptive set theory, but without using effective descriptive set theory. He replaces it with some versions of Theorem 1.2. In particular, he can prove Theorem 1.2 without effective descriptive set theory.\bigskip

  When $A$ is Borel, it is natural to ask about the relation between the Borel class of $A$ and that of the coloring $f$ when Theorem 1.2.(a) holds. This leads to consider $\borxi$-measurable countable colorings (or equivalently $\boraxi$-measurable countable colorings). We have the following conjecture:\bigskip
   
\noindent\bf Conjecture 1\it\ Let $1\!\leq\!\xi\! <\!\omega_1$. Then there are\smallskip

\noindent - a $0$-dimensional Polish space $\mathbb{X}_\xi$,\smallskip

\noindent - an analytic relation $\mathbb{A}_{\xi}$ on $\mathbb{X}_\xi$\smallskip

\noindent such that for any $0$-dimensional Polish space $X$, and for any analytic relation $A$ on $X$, exactly one of the following holds:\smallskip  

\noindent (a) $(X,A)\preceq_{\borxi}\big(\omega ,\neg\Delta (\omega )\big)$,\smallskip  

\noindent (b) $(\mathbb{X}_\xi ,\mathbb{A}_{\xi})\preceq_{\boraone}(X,A)$.\rm\bigskip

 We will prove it when $1\!\leq\!\xi\!\leq\! 2$, and in these cases we do not have to assume that $A$ is analytic. A sequence $s\!\in\! 3^{<\omega}$ will be said to be $suitable$ if $s\! =\!\emptyset$ or 
$s(\vert s\vert\! -\! 1)\! =\! 2$. We will have $\mathbb{X}_2\! :=\! 3^\omega$ and $\mathbb{A}_2\! :=\!
\big\{ (s0\alpha ,s1\beta )\mid s\mbox{ suitable}\wedge\alpha ,\beta\!\in\! 2^\omega\big\}$.

\vfill\eject 

 We saw that the study of the Borel countable colorings is highly related to the study of countable unions of Borel rectangles, and gave some motivation for studying $\boraxi$-measurable countable colorings. This motivates the study of countable unions of $\boraxi$ rectangles. Another motivation is that 
$(X,A)\preceq_{\borxi}\big(\omega ,\neg\Delta (\omega )\big)$ is equivalent to the fact that $\Delta (X)$ can be separated from $A$ by a $(\boraxi\!\times\!\boraxi )_\sigma$ set, by the generalized reduction property for the class $\boraxi$ (see 22.16 in [K]).\bigskip
   
\noindent\bf Conjecture 2\it\ Let $1\!\leq\!\xi\! <\!\omega_1$. Then there are\smallskip

\noindent - $0$-dimensional Polish spaces $\mathbb{X}^0_\xi ,\mathbb{X}^1_\xi$,\smallskip

\noindent - disjoint analytic subsets $\mathbb{A}^0_{\xi},\mathbb{A}^1_{\xi}$ of 
$\mathbb{X}^0_\xi\!\times\!\mathbb{X}^1_\xi$\smallskip

\noindent such that for any Polish spaces $X,Y$, and for any pair $A_0,A_1$ of disjoint analytic subsets of $X\!\times\! Y$, exactly one of the following holds:\smallskip  

\noindent (a) The set $A_0$ can be separated from $A_1$ by a $(\boraxi\!\times\!\boraxi )_\sigma$ set,\smallskip  

\noindent (b) $(\mathbb{X}^0_\xi ,\mathbb{X}^1_\xi ,\mathbb{A}^0_\xi ,\mathbb{A}^1_\xi )\leq 
(X,Y,A_0,A_1)$.\rm\bigskip

 It is easy to prove this when $\xi\! =\! 1$. Our main result is that Conjecture 2 holds when 
$\xi\! =\! 2$. We now describe our minimum example 
$(\mathbb{X}^0_2,\mathbb{X}^1_2,\mathbb{A}^0_2,\mathbb{A}^1_2)$.\bigskip

\noindent\bf Notation.\rm\ We put $\mathbb{X}^0_2\! :=\! 3^\omega\!\setminus\!\big\{ s1\beta\mid 
s\mbox{ suitable}\wedge\beta\!\in\! 2^\omega\big\}$, 
$\mathbb{X}^1_2\! :=\! 3^\omega\!\setminus\!\big\{ s0\alpha\mid s\mbox{ suitable}\wedge
\alpha\!\in\! 2^\omega\big\}$,  
$\mathbb{A}^0_2\! :=\!\Delta (\mathbb{X}^0_2\cap\mathbb{X}^1_2)$ and 
$\mathbb{A}^1_2\! :=\!\mathbb{A}_2\! :=\!\big\{ (s0\alpha ,s1\beta )\mid s\mbox{ suitable}\wedge
\alpha ,\beta\!\in\! 2^\omega\big\}$.\bigskip

 We use effective descriptive set theory, and give effective strengthenings of our results. The reader should see [M] for basic notions of effective descriptive set theory. In particular, we will see that to test whether an analytic relation has a $\boraxi$-measurable countable coloring, it is enough to test countably many partitions instead of continuum many. We will use the topology $T_2$ generated by the 
$\Ana\cap\bormone$ subsets of a recursively presented Polish space (introduced in [Lo1]). Our main result can be strengthened as follows (see [L3]).

\begin{thm} Let $X,Y$ be recursively presented Polish spaces, and 
$A_0,A_1$ disjoint $\Ana$ subsets of $X\!\times\! Y$. The following are equivalent:\smallskip

\noindent (a) The set $A_0$ cannot be separated from $A_1$ by a 
$(\boratwo\!\times\!\boratwo )_\sigma$ set.\smallskip

\noindent (b) The set $A_0$ cannot be separated from $A_1$ by a 
$\Borel\cap (\boratwo\!\times\!\boratwo )_\sigma$ set.\smallskip

\noindent (c) The set $A_0$ cannot be separated from $A_1$ by a 
$\boraone (T_2\!\times\! T_2)$ set.\smallskip

\noindent (d) $A_0\cap\overline{A_1}^{T_2\times T_2}\!\not=\!\emptyset$.\smallskip

\noindent (e) $(\mathbb{X}^0_2,\mathbb{X}^1_2,\mathbb{A}^0_2,\mathbb{A}^1_2)\leq 
(X,Y,A_0,A_1)$.\end{thm}

\vfill\eject

\section{$\!\!\!\!\!\!$ Some general effective facts}\indent

 One can hope for an effective strengthening of Conjecture 1:
    
\noindent\bf Effective conjecture 1\it\ Let $1\!\leq\!\xi\! <\!\omega_1$. Then there are\smallskip

\noindent - a $0$-dimensional Polish space $\mathbb{X}_\xi$,\smallskip

\noindent - an analytic relation $\mathbb{A}_{\xi}$ on $\mathbb{X}_\xi$\smallskip

\noindent such that 
$(\mathbb{X}_\xi ,\mathbb{A}_{\xi})\not\preceq_{\borxi}\big(\omega ,\neg\Delta (\omega )\big)$ and 
for any $\alpha\!\in\!\omega^\omega$ with $1\!\leq\!\xi\! <\!\omega_1^\alpha$, for any 
$0$-dimensional recursively in $\alpha$ presented Polish space $X$, and for any $\Ana (\alpha )$ relation $A$ on $X$, one of the following holds:\smallskip  

\noindent (a) $(X,A)\preceq_{\Borel (\alpha )\cap\borxi}
\big(\omega ,\neg\Delta (\omega )\big)$,\smallskip  

\noindent (b) $(\mathbb{X}_\xi ,\mathbb{A}_{\xi})\preceq_{\boraone}(X,A)$.\rm\bigskip

 We will see that this effective conjecture is true when $1\!\leq\!\xi\!\leq\! 2$. The following statement is a corollary of this effective conjecture, and is in fact a theorem: 
   
\begin{thm} Let $1\!\leq\!\xi\! <\!\omega_1^{\mbox{CK}}$, $X$ a $0$-dimensional recursively presented Polish space, and $A$ a $\Ana$ relation on $X$. We assume that 
$(X,A)\preceq_{\borxi}\big(\omega ,\neg\Delta (\omega )\big)$. Then 
$(X,A)\preceq_{\Borel\cap\borxi}\big(\omega ,\neg\Delta (\omega )\big)$.\end{thm}

 A consequence of this is that to test whether an analytic relation has a $\boraxi$-measurable countable coloring, it is enough to test countably many partitions instead of continuum many. Another consequence is the equivalence between Conjecture 1 and the Effective conjecture 1. We have in fact preliminary results that will help us to prove also the equivalence between (a)-(d) in Theorem 1.5, in the general case. 
 
\begin{lem} Let $1\!\leq\!\xi\! <\!\omega_1^{\mbox{CK}}$, $X, Y$ recursively presented Polish spaces, and $A\!\in\!\Ana (X)\cap\boraxi$, $B\!\in\!\Ana (Y)\cap\boraxi$ and $C\!\in\!\Ana (X\!\times\! Y)$ disjoint from $A\!\times\! B$. Then there are $A',B'\!\in\!\Borel\cap\boraxi$ such that $A'\!\times\! B'$ separates 
$A\!\times\! B$ from $C$.\end{lem}
 
\noindent\bf Proof.\rm\ Note that $A$ and $\{ x\!\in\! X\mid\exists y\!\in\! B~~(x,y)\!\in\! C\}$ are disjoint 
$\Ana$ sets, separable by a $\boraxi$ subset of $X$. By Theorems 1.A and 1.B in [Lo1], there is 
$A'\!\in\!\Borel\cap\boraxi$ separating these two sets. Similarly, $B$ and 
$\{ y\!\in\! Y\mid\exists x\!\in\! A'~~(x,y)\!\in\! C\}$ are disjoint $\Ana$ sets, and there is 
$B'\!\in\!\Borel\cap\boraxi$ separating these two sets.\hfill{$\square$}

\begin{thm} Let $1\!\leq\!\xi\! <\!\omega_1^{\mbox{CK}}$, $X, Y$ recursively presented Polish spaces, and $A_0,A_1$ disjoint $\Ana$ subsets of $X\!\times\! Y$. We assume that $A_0$ is separable from $A_1$ by a $\big(\boraxi\!\times\!\boraxi\big)_\sigma$ set. Then $A_0$ is separable from $A_1$ by a 
$\Borel\cap\big( (\Borel\cap\boraxi )\!\times\! (\Borel\cap\boraxi)\big)_\sigma$ set.\end{thm}

\noindent\bf Proof.\rm\ By Example 2 of Chapter 3 in [Lo2], the family $\big( N(n,X)\big)_{n\in\omega}$ is regular without parameter. By Corollary 2.10 in [Lo2], $\bormxi (X)$, as well as 
$\boraxi (X)\! =\!\big(\bigcup_{\eta <\xi}~\borme (X)\big)_\sigma$,  are regular without parameter. By Theorem 2.12 in [Lo2], $\boraxi (X)\!\times\!\boraxi (Y)$ is also regular without parameter. By Theorem 2.8 in [Lo2], the family $\Phi\! :=\!\big(\boraxi (X)\!\times\!\boraxi (Y)\big)_\sigma$ is separating which imply the existence of $S\!\in\!\Borel\cap\Phi$ separating $A_0$ from $A_1$.

\vfill\eject

 With the notation of [Lo2], let $n$ be an integer with $(0^\infty ,n)\!\in\! W$ and $C_{0^\infty ,n}\! =\! S$. Then $(0^\infty ,n)$ is in $W_\Phi$, which by Theorem 2.8.(ii) in [Lo2] is 
$$\left\{ (\alpha ,n)\!\in\! W\mid\exists\beta\!\in\!\Borel (\alpha )~~\forall m\!\in\!\omega ~~
\big(\alpha ,\beta (m)\big)\!\in\! W_{\boraxi (X)\times\boraxi (Y)}\wedge
C_{\alpha ,n}\! =\!\bigcup_{m\in\omega}~C_{\alpha ,\beta (m)}\right\}.$$
This implies that $S\!\in\!\Borel\cap\big(\Borel\cap (\boraxi\!\times\!\boraxi )\big)_\sigma$. It remains to check that 
$\Borel\cap (\boraxi\!\times\!\boraxi )\! =\! (\Borel\cap\boraxi )\!\times\! (\Borel\cap\boraxi )$. The second set is clearly a subset of the first one. So assume that 
$R\! =\! A\!\times\! B\!\in\!\Borel\cap (\boraxi\!\times\!\boraxi )$. We may assume that $R$ is not empty. Then the projections $A$, $B$ are $\Ana$ since $R\!\in\!\Borel$. Lemma 2.2 gives 
$A',B'\!\in\!\Borel\cap\boraxi$ with $A\!\times\! B\!\subseteq\! A'\!\times\! B'\!\subseteq\! R\! =\! A\!\times\! B$. \hfill{$\square$}\bigskip

 Recall that if $A$ is a relation on $X$ and $D\!\subseteq\! X$, then $D$ is $A$-$discrete$ if 
$A\cap D^2\! =\!\emptyset$.\bigskip

\noindent\bf Proof of Theorem 2.1.\rm\ We apply Theorem 2.3 to $Y\! :=\! X$, $A_0\! :=\!\Delta (X)$ and $A_1\! :=\! A$. As 
$$(X,A)\preceq_{\borxi}\big(\omega ,\neg\Delta (\omega )\big)\mbox{,}$$ 
$\Delta (X)$ is separable from $A$ by a $(\boraxi\!\times\!\boraxi )_\sigma$ set. Theorem 2.3 gives 
$C_n,D_n\!\in\!\Borel\cap\boraxi$ such that $S\! :=\!\bigcup_{n\in\omega}~C_n\!\times\! D_n\!\in\!\Borel$ separates $\Delta (X)$ from $A$. As the set of codes for $\Borel\cap\boraxi$ subsets of $X$ is $\Ca$ (see Proposition 1.4 in [Lo1]), the $\Borel$-selection theorem and the separation theorem imply that we may assume that the sequences $(C_n)$ and $(D_n)$ are $\Borel$. Note that $(C_n\cap D_n)$ is a $\Borel$ covering of $X$ into $A$-discrete $\Borel\cap\boraxi$ sets. As $X$ is $0$-dimensional we can reduce this covering into a $\Borel$ covering $(\Delta_n)$ of $X$ into $\Borel\cap\boraxi$ sets, which are in fact $\borxi$. This gives the desired partition.\hfill{$\square$}\bigskip

\noindent\bf Notation.\rm\ Following [Lo1], we define the following topologies on a $0$-dimensional recursively in $\alpha$ presented Polish space $X$, for any $\alpha\!\in\!\omega^\omega$. Let 
$T_1 (\alpha )$ be the usual topology on $X$, and for $2\!\leq\!\xi\! <\!\omega_1$, $T_\xi (\alpha )$ be the topology generated by the $\Ana (\alpha )\cap\bormlxi$ subsets of $X$. The next proposition gives a reformulation of the inequality 
$(X,A)\preceq_{\Borel (\alpha )\cap\borxi}\big(\omega ,\neg\Delta (\omega )\big)$ of the Effective 
conjecture 1.

\begin{prop} Let $1\!\leq\!\xi\! <\!\omega_1^{\mbox{CK}}$, $X$ a $0$-dimensional recursively presented Polish space, and $A$ a $\Ana$ relation on $X$. Then 
$(X,A)\preceq_{\Borel\cap\borxi}\big(\omega ,\neg\Delta (\omega )\big)$ is equivalent to 
$\Delta (X)\cap\overline{A}^{T_\xi\times T_\xi}\! =\!\emptyset$.\end{prop}

\noindent\bf Proof.\rm\ Assume first that 
$(X,A)\preceq_{\Borel\cap\borxi}\big(\omega ,\neg\Delta (\omega )\big)$. Then there is a partition $(B_n)$ of $X$ into $A$-discrete $\Borel\cap\borxi$ sets. In particular, Theorem 1.A in [Lo1] implies that  $B_n$ is a countable union of $\Borel\cap\bormlxi$ sets if $\xi\!\geq\! 2$. In particular, $B_n$ is $T_\xi$-open and $\Delta (X)$ is disjoint from $\overline{A}^{T_\xi\times T_\xi}$.\bigskip

 Conversely, assume that $\Delta (X)\cap\overline{A}^{T_\xi\times T_\xi}\! =\!\emptyset$. Then each element $x$ of $X$ is contained in a $A$-discrete $\Ana\cap\bormlxi$ set (basic clopen set if 
$\xi\! =\! 1$). Lemma 2.2 implies that each element $x$ of $X$ is in fact contained in a $A$-discrete $\Borel\cap\bormlxi$ set if $\xi\!\geq\! 2$. It remains to apply Proposition 1.4 in [Lo1] and the $\Borel$-selection theorem to get the desired partition.\hfill{$\square$}

\vfill\eject

 One can also hope for an effective strengthening of Conjecture 2 generalizing Theorem 1.5:\bigskip

\noindent\bf Effective conjecture 2\it\ Let $1\!\leq\!\xi\! <\!\omega_1$. Then there are\smallskip

\noindent - $0$-dimensional Polish spaces $\mathbb{X}^0_\xi ,\mathbb{X}^1_\xi$,\smallskip

\noindent - disjoint analytic subsets $\mathbb{A}^0_{\xi},\mathbb{A}^1_{\xi}$ of the space 
$\mathbb{X}^0_\xi\!\times\!\mathbb{X}^1_\xi$, not separable by a $(\boraxi\!\times\boraxi )_\sigma$ set,\smallskip

\noindent such that for any $\alpha\!\in\!\omega^\omega$ such that $1\!\leq\!\xi\! <\!\omega_1^\alpha$, for any recursively in $\alpha$ presented Polish spaces $X,Y$, and for any pair 
$A_0,A_1$ of disjoint $\Ana (\alpha )$ subsets of $X\!\times\! Y$, the following are equivalent:\smallskip

\noindent (a) The set $A_0$ cannot be separated from $A_1$ by a 
$(\boraxi\!\times\!\boraxi )_\sigma$ set.\smallskip

\noindent (b) The set $A_0$ cannot be separated from $A_1$ by a 
$\Borel (\alpha )\cap (\boraxi\!\times\!\boraxi )_\sigma$ set.\smallskip

\noindent (c) The set $A_0$ cannot be separated from $A_1$ by a 
$\boraone\big( T_\xi (\alpha )\!\times\! T_\xi (\alpha )\big)$ set.\smallskip

\noindent (d) $A_0\cap\overline{A_1}^{T_\xi (\alpha )\times T_\xi (\alpha )}\!\not=\!\emptyset$.\smallskip

\noindent (e) $(\mathbb{X}^0_\xi ,\mathbb{X}^1_\xi ,\mathbb{A}^0_\xi ,\mathbb{A}^1_\xi )\leq 
(X,Y,A_0,A_1)$.\rm\bigskip

 In fact, the statements (a)-(d) are indeed equivalent:

\begin{thm} Let $1\!\leq\!\xi\! <\!\omega_1^{\mbox{CK}}$, $X, Y$ recursively presented Polish spaces, and $A_0,A_1$ disjoint $\Ana$ subsets of $X\!\times\! Y$. The following are equivalent:\smallskip

\noindent (a) The set $A_0$ cannot be separated from $A_1$ by a 
$(\boraxi\!\times\!\boraxi )_\sigma$ set.\smallskip

\noindent (b) The set $A_0$ cannot be separated from $A_1$ by a 
$\Borel\cap (\boraxi\!\times\!\boraxi )_\sigma$ set.\smallskip

\noindent (c) The set $A_0$ cannot be separated from $A_1$ by a 
$\boraone (T_\xi\!\times\! T_\xi )$ set.\smallskip

\noindent (d) $A_0\cap\overline{A_1}^{T_\xi\times T_\xi}\!\not=\!\emptyset$.\end{thm}

\noindent\bf Proof.\rm\ Theorem 2.3 implies that (a) is indeed equivalent to (b). It also implies, using the proof of Proposition 2.4, that (c) implies (a), and the converse is clear. It is also clear that (c) and (d) are equivalent.\hfill{$\square$}\bigskip

 A consequence of this is that Conjecture 2 and the Effective conjecture 2 are equivalent.

\section{$\!\!\!\!\!\!$ The case $\xi\! =\! 1$}\indent

 We set $\mathbb{X}_1\! :=\! 2^\omega$ and $\mathbb{A}_1\! :=\!\{ (0^{2k+1}1\alpha ,0^{2k}1\beta )\mid k\!\in\!\omega\wedge\alpha ,\beta\!\in\! 2^\omega\}$.

\begin{lem} The space $\mathbb{X}_1$ is a $0$-dimensional metrizable compact space, $\mathbb{A}_1$ is a $\boraone$ relation on $\mathbb{X}_1$, and $(\mathbb{X}_1,\mathbb{A}_1)\not\preceq_{\borone}\big(\omega ,\neg\Delta (\omega )\big)$.\end{lem}
 
\noindent\bf Proof.\rm\ The first two assertions are clear. We argue by contradiction for the last assertion, which gives $f\! :\!\mathbb{X}_1\!\rightarrow\!\omega$ continuous with 
$f(\alpha )\!\not=\! f(\beta )$ if $(\alpha ,\beta )\!\in\!\mathbb{A}_1$. We set 
$C_n\! :=\! f^{-1}(\{ n\})$, so that $(C_n)_{n\in\omega}$ is a partition of $\mathbb{X}_1$ into $\mathbb{A}_1$-discrete $\borone$ sets. Choose $n$ with $0^\infty\!\in\! C_n$. Then 
$0^i\alpha\!\in\! C_n$ if $i$ is big enough. This gives an integer $k$ with 
$0^{2k+1}1^\infty ,0^{2k}1^\infty\!\in\! C_n$, and 
$(0^{2k+1}1^\infty ,0^{2k}1^\infty )\!\in\!\mathbb{A}_1\cap C_n^2$, which is absurd.\hfill{$\square$}

\vfill\eject

\begin{thm} Let $X$ be a $0$-dimensional Polish space, and $A$ a relation on $X$. Then exactly one of the following holds:\smallskip  

\noindent (a) $(X,A)\preceq_{\borone}\big(\omega ,\neg\Delta (\omega )\big)$,\smallskip  

\noindent (b) $(\mathbb{X}_1,\mathbb{A}_1)\preceq_{\boraone}(X,A)$.\smallskip

 Moreover, this is not true, even if $A$ is analytic, if $X$ is not $0$-dimensional, and we cannot have $f$ one-to-one in (b) (with this couple $(\mathbb{X}_1,\mathbb{A}_1)$ or any other).\end{thm}
 
\noindent\bf Proof.\rm\ Note first that (a) and (b) cannot hold simultaneously, by Lemma 3.1. We enumerate a basis $\big( N(n,X)\big)_{n\in\omega}$ for the topology of $X$ made of clopen sets. Assume that (a) does not hold. We build\medskip
 
 \noindent - an increasing sequence of integers $(n_k)_{k\in\omega}$,\smallskip
 
 \noindent - a sequence $(x_p)_{p\in\omega}$ of points of $X$.\bigskip
 
 We want these objects to satisfy the following conditions:
$$\begin{array}{ll}
& (1)~(x_{2k},x_{2k+1})\!\in\! A\cap N(n_k,X)^2\cr
& (2)~N(n_{k+1},X)\!\subseteq\! N(n_k,X)\cr
& (3)~\mbox{diam}\big( N(n_k,X)\big)\!\leq\! 2^{-k}\cr
& (4)~\mbox{There is no covering of }N(n_k,X)
\mbox{ into }A\mbox{-discrete clopen subsets of }X
\end{array}$$ 
$\bullet$ Assume that this is done. Then we can define a point $x$ of $X$ by 
$\{ x\}\! =\!\bigcap_{k\in\omega}~N(n_k,X)$. Note that $(x_p)$ tends to $x$. We define 
$f\! :\!\mathbb{X}_1\!\rightarrow\! X$ by $f(0^\infty )\! :=\! x$, $f(0^{2k+1}1\alpha )\! :=\! x_{2k}$ and 
$f(0^{2k}1\beta )\! :=\! x_{2k+1}$. Note that $f$ is continuous. Moreover, 
$\big( f(0^{2k+1}1\alpha ),f(0^{2k}1\beta )\big)\! =\! (x_{2k},x_{2k+1})\!\in\! A$, so that (b) holds.\bigskip

\noindent $\bullet$ Let us prove that the construction is possible. We set $N(n_{-1},X)\! :=\! X$. Assume that $(n_k)_{k<l}$ and $(x_{2k},x_{2k+1})_{k<l}$ satisfying (1)-(4) have been constructed, which is the case for $l\! =\! 0$. We choose a covering of $N(n_{l-1},X)$ with basic clopen sets of diameter at most $2^{-l}$, contained in $N(n_{l-1},X)$. Then one of these basic sets, say $N(n_l,X)$, satisfies (4). It remains to choose $(x_{2l},x_{2l+1})$ in the set $A\cap N(n_l,X)^2$.\bigskip

\noindent $\bullet$ Consider now $X\! :=\!\mathbb{R}$ and $A\! :=\!\{ (0,1)\}$. Then (a) does not hold since $\mathbb{R}$ is connected. If (b) holds, then we must have $f(0^{2k+1}1\alpha )\! =\! 0$ and 
$f(0^{2k}1\beta )\! =\! 1$. By continuity of $f$, we get $f(0^\infty )\! =\! 0\! =\! 1$.\bigskip

 This would be the same with any $(\mathbb{X}_1,\mathbb{A}_1)$. Indeed, as 
$(\mathbb{X}_1,\mathbb{A}_1)\not\preceq_{\borone}\big(\omega ,\neg\Delta (\omega )\big)$, we have 
$\overline{\Pi_0[\mathbb{A}_1]}\cap\overline{\Pi_1[\mathbb{A}_1]}\!\not=\!\emptyset$, since otherwise there would be a clopen subset $C$ of $\mathbb{X}_1$ separating $\overline{\Pi_0[\mathbb{A}_1]}$ from $\overline{\Pi_1[\mathbb{A}_1]}$, and we would have 
$\Delta (\mathbb{X}_1)\!\subseteq\! C^2\cup (\neg C)^2\!\subseteq\!\neg\mathbb{A}_1$. So we can choose $x\!\in\!\overline{\Pi_0[\mathbb{A}_1]}\cap\overline{\Pi_1[\mathbb{A}_1]}$, 
$x_{2k}\!\in\!\Pi_0[\mathbb{A}_1]$ such that $(x_{2k})$ tends to $x$, 
$y_{2k+1}\!\in\!\Pi_1[\mathbb{A}_1]$ such that $(y_{2k+1})$ tends to $x$, 
$y_{2k}$ with $(x_{2k},y_{2k})\!\in\!\mathbb{A}_1$, and 
$x_{2k+1}$ with $(x_{2k+1},y_{2k+1})\!\in\!\mathbb{A}_1$. Then $f(x_{2k})\! =\! 0$, $f(y_{2k+1})\! =\! 1$ and we conclude as before.\bigskip

\noindent $\bullet$ Consider $X\! :=\! 2^\omega$ and 
$A\! :=\!\{ 0^\infty\}\!\times\! (2^\omega\!\setminus\!\{ 0^\infty\} )$. Then (a) does not hold since if a clopen subset $C$ of $2^\omega$ contains $0^\infty$, then it contains also $\alpha\!\not=\! 0^\infty$, so that 
$(0^\infty ,\alpha )\!\in\! A\cap C^2$. If (b) holds, then $f(0^{2k+1}1\alpha )\! =\! 0^\infty$ for each integer $k$ and $f$ is not one-to-one.\bigskip

 This argument works as soon as $\Pi_0[\mathbb{A}_1]$ has at least two elements. If we argue in the other factor, then we see that an example $(\mathbb{X}_1,\mathbb{A}_1)$ with injectivity must satisfy that $\mathbb{A}_1$ is a singleton $\{ (\alpha ,\beta )\}$. As 
$(\mathbb{X}_1,\mathbb{A}_1)\preceq_{\boraone}(2^\omega ,\mathbb{G}_0)$, $\alpha\!\not=\!\beta$. So take a clopen subset $C$ of $\mathbb{X}_1$ containing $\alpha$ but not $\beta$. Then 
$\Delta (\mathbb{X}_1)\!\subseteq\! C^2\cup (\neg C)^2\!\subseteq\!\neg\mathbb{A}_1$.\hfill{$\square$}

\vfill\eject

\begin{prop} Conjecture 2 holds for $\xi\! =\! 1$.\end{prop}

\noindent\bf Proof.\rm\ We set $\mathbb{X}^\varepsilon_1\! :=\!\mathbb{X}_1$, 
$\mathbb{A}^0_1\! :=\!\{ (0^\infty ,0^\infty )\}$ and $\mathbb{A}^1_1\! :=\!\mathbb{A}_1$. If 
$(x,y)\!\in\! A_0\cap\overline{A_1}$, then choose $(x_k,y_k)$ in $A_1$ tending to $(x,y)$, and set 
$f(0^\infty )\! :=\! x$, $g(0^\infty )\! :=\! y$, $f(0^{2k+1}1\alpha )\! :=\! f(0^{2k}1\beta )\! :=\! x_k$, 
$g(0^{2k+1}1\alpha )\! :=\! g(0^{2k}1\beta )\! :=\! y_k$.\hfill{$\square$}

\section{$\!\!\!\!\!\!$ The case $\xi\! =\! 2$}

\begin{lem} The space $\mathbb{X}_2$ is a $0$-dimensional metrizable compact space, $\mathbb{A}_2$ is a $\boratwo$ relation on $\mathbb{X}_2$, and $(\mathbb{X}_2,\mathbb{A}_2)\not\preceq_{\bortwo}\big(\omega ,\neg\Delta (\omega )\big)$.\end{lem}
 
\noindent\bf Proof.\rm\ The first two assertions are clear. We argue by contradiction for the last assertion, which gives $f\! :\!\mathbb{X}_2\!\rightarrow\!\omega$ $\bortwo$-measurable with 
$f(\alpha )\!\not=\! f(\beta )$ if $(\alpha ,\beta )\!\in\!\mathbb{A}_2$. We set 
$C_n\! :=\! f^{-1}(\{ n\})$, so that $(C_n)_{n\in\omega}$ is a partition of $\mathbb{X}_2$ into 
$\mathbb{A}_2$-discrete $\bortwo$ sets. By Baire's theorem, there are an integer $n$ and 
$s\!\in\! 2^{<\omega}$ such that $C_n$ contains the basic clopen set $N_s$. Then 
$(s20^\infty ,s21^\infty )\!\in\!\mathbb{A}_2\cap C_n^2$, which is absurd.\hfill{$\square$}\bigskip

 We have a stronger result than Conjecture 1, in the sense that we do not need any regularity assumption on $A$, neither that $X$ is $0$-dimensional:

\begin{thm} (Lecomte-Zelen\'y) Let $X$ be a Polish space, and $A$ a relation on $X$. Then exactly one of the following holds:\smallskip  

\noindent (a) $(X,A)\preceq_{\bortwo}\big(\omega ,\neg\Delta (\omega )\big)$,\smallskip  

\noindent (b) $(\mathbb{X}_2,\mathbb{A}_2)\preceq_{\boraone}(X,A)$.\end{thm}
 
\noindent\bf Proof.\rm\ Note first that (a) and (b) cannot hold simultaneously, by Lemma 4.1. If $A$ is not a digraph, then choose $x$ with $(x,x)\!\in\! A$, and put $f(\alpha )\! :=\! x$. So we may assume that $A$ is a digraph. We set 
$$U\! :=\!\bigcup\Big\{ V\!\in\!\boraone (X)\mid\exists D\!\in\!\boratwo (\omega\!\times\! X)~~
V\!\subseteq\!\bigcup_{p\in\omega}~D_p\wedge\forall p\!\in\!\omega ~~
A\cap (D_p\!\times\! D_p)\! =\!\emptyset\Big\}.$$
\bf Case 1.\rm\ $U\! =\! X$.\bigskip

 There is a countable covering of $X$ into $A$-discrete $\boratwo$ sets. We just have to reduce them to get a partition showing that (a) holds.\bigskip

\noindent\bf Case 2.\rm\ $U\!\not=\! X$.\bigskip

 Then $Y\! :=\! X\!\setminus\! U$ is a nonempty closed subset of $X$.\bigskip

\noindent\bf Claim\it\ If $\emptyset\!\not=\! W\!\in\!\boraone (Y)$, then there is no $\boratwo$ subset of 
$\omega\!\times\! X$ whose sections are $A$-discrete and cover $W$. In particular, $W$ is not $A$-discrete.\rm\bigskip
 
 We argue by contradiction. Let $y\!\in\! W$, and $Z$ an open subset of $X$ with $Z\cap Y\! =\! W$. As 
$Z\cap U$ can be covered with some $\bigcup_{p\in\omega}~D_p$'s, so is $Z$. Thus 
$Z\!\subseteq\! U$, so that $y\!\in\! Z\cap Y\!\subseteq\! U\!\setminus\! U\! =\!\emptyset$, which is the desired contradiction.\hfill{$\diamond$}\bigskip

 We construct a sequence $(V_s)_{s\in 3^{<\omega}}$ of open subsets of $Y$, and a sequence $(x_s)_{s\in 3^{<\omega}}$ of points of $Y$. We want these objects to satisfy the following conditions:
$$\begin{array}{ll}
& (1)~x_s\!\in\! V_s\cr
& (2)~\overline{V_{s\varepsilon}}\!\subseteq\! V_s\cr
& (3)~\mbox{diam}(V_s)\!\leq\! 2^{-\vert s\vert}\cr
& (4)~(x_{s0},x_{s1})\!\in\! A\mbox{ if }s\mbox{ is suitable}\cr
& (5)~x_{s\varepsilon}\! =\! x_s\mbox{ if }\varepsilon\! =\! 2\vee s\mbox{ is not suitable}
\end{array}$$ 
$\bullet$ Assume that this is done. We define $f\! :\! 3^\omega\!\rightarrow\! Y\!\subseteq\! X$ by 
$\{ f(\alpha )\}\! :=\!\bigcap_{k\in\omega}~\overline{V_{\alpha\vert k}}\! =\!
\bigcap_{k\in\omega}~V_{\alpha\vert k}$, so that $f$ is continuous. Note that 
$f(\alpha )$ is the limit of $x_{\alpha\vert k}$, and that
$$x_{s\varepsilon}\! =\! x_{s\varepsilon (\alpha\vert 1)}\! =\! ...\! =\! x_{s\varepsilon (\alpha\vert (q+1))}$$
for each $(s,\varepsilon ,\alpha )\!\in\! 3^{<\omega}\!\times\! 2\!\times\! 2^\omega$. Thus 
$f(s\varepsilon\alpha )\! =\!\mbox{lim}_{q\rightarrow\infty}~x_{s\varepsilon (\alpha\vert q)}
\! =\! x_{s\varepsilon}$ and 
$$\big( f(s0\alpha ),f(s1\beta )\big)\! =\! (x_{s0},x_{s1})\!\in\! A.$$ 
So (b) holds.\bigskip

\noindent $\bullet$ Let us prove that the construction is possible. We choose $x_\emptyset\!\in\! Y$ and 
an open neighborhood $V_\emptyset$ of $x_\emptyset$ in $Y$, of diameter at most $1$. Assume that 
$(V_s)_{s\in 3^{\leq l}}$ and $(x_s)_{s\in 3^{\leq l}}$ satisfying (1)-(5) have been constructed, which is the case for $l\! =\! 0$.\bigskip

 An application of the Claim gives $(x_{s0},x_{s1})\!\in\! A\cap V_s^2$ if $s$ is suitable. We satisfy (5), so that the definition of the $x_s$'s is complete. Note that  $x_s\!\in\! V_{s\vert l}$ if $s\!\in\! 3^{l+1}$.\bigskip

 We choose an open neighborhood $V_s$ of $x_s$ in $Y$, of diameter at most $2^{-l-1}$, ensuring the inclusion $\overline{V_s}\!\subseteq\! V_{s\vert l}$. This finishes the proof.\hfill{$\square$}\bigskip

\noindent\bf Remark.\rm\ We cannot replace $(\mathbb{X}_2,\mathbb{A}_2)$ with $\big( 2^\omega ,
\big\{ (s0\alpha ,s1\beta )\mid s\!\in\! 2^{<\omega}\wedge\alpha ,\beta\!\in\! 2^\omega\big\}\big)$. Indeed, otherwise we get $f\! :\! 2^\omega\!\rightarrow\! 3^\omega$ continuous with 
$$\big\{ (s0\alpha ,s1\beta )\mid s\!\in\! 2^{<\omega}\wedge\alpha ,\beta\!\in\! 2^\omega\big\}
\!\subseteq\! (f\!\times\!f)^{-1}(\big\{ (s0\alpha ,s1\beta )\mid
s\mbox{ suitable }\wedge\alpha ,\beta\!\in\! 2^\omega\big\}).$$ 
Thus $\big( f(0^\infty ),f(0^k1^\infty )\big)\! =\! (s_k0\alpha_k,s_k1\beta_k)\! =\! (s_00\alpha_0,s_01\beta_k)$. But $f(0^\infty )\! =\! s_00\alpha_0$ is the limit of 
$f(0^k1^\infty )\! =\! s_01\beta_k$, which cannot be. This shows that it is useful to take $3$ instead of $2$.\bigskip

 Now we come to the proof of our main theorem.

\begin{lem} The spaces $\mathbb{X}^0_2,\mathbb{X}^1_2$ are $0$-dimensional Polish spaces, 
$\mathbb{A}^0_2,\mathbb{A}^1_2$ are disjoint analytic subsets of 
$\mathbb{X}^0_2\!\times\!\mathbb{X}^1_2$, and are not separable by a 
$(\boratwo\!\times\!\boratwo )_\sigma$ set.\end{lem}
 
\noindent\bf Proof.\rm\ The first two assertions are clear since $\mathbb{X}^0_2,\mathbb{X}^1_2$ are 
$G_\delta$ subsets of $3^\omega$, $\mathbb{A}^0_2,\mathbb{A}^1_2$ have disjoint projections, 
$\mathbb{A}^0_2\! =\!\Delta (3^\omega )\cap (\mathbb{X}^0_2\!\times\!\mathbb{X}^1_2)$ is closed and 
$\mathbb{A}^1_2$ is $\boratwo$. We argue by contradiction for the last assertion, which gives 
$C_n\!\in\!\bormone (\mathbb{X}^0_2)$ and $D_n\!\in\!\bormone (\mathbb{X}^1_2)$ with 
$\mathbb{A}^0_2\!\subseteq\!\bigcup_{n\in\omega}~(C_n\!\times\! D_n)\!\subseteq\!
\neg\mathbb{A}^1_2$. In particular, $\mathbb{X}^0_2\cap\mathbb{X}^1_2\! =\!\bigcup_{n\in\omega}~C_n\cap D_n$, and Baire's theorem gives $n$ and $s\!\in\! 3^{<\omega}$ such that the inclusion $N_s\cap\mathbb{X}^0_2\cap\mathbb{X}^1_2\!\subseteq\! C_n\cap D_n$ holds. Note that $N_s\cap\mathbb{X}^0_2\!\subseteq\! C_n$ and $N_s\cap\mathbb{X}^1_2\!\subseteq\! D_n$. Then $(s20^\infty ,s21^\infty)\!\in\!\mathbb{A}^1_2\cap (C_n\!\times\! D_n)$, which is absurd.\hfill{$\square$}

\vfill\eject

\noindent\bf Remark.\rm\ This proof shows that the spaces $\mathbb{X}^0_2,\mathbb{X}^1_2$ of Conjecture 2 cannot be both compact, which is quite unusual in this kind of dichotomy (even if it was already the case in [L2]). Indeed, our example shows that $\mathbb{A}^0_2,\mathbb{A}^1_2$ must be separable by a closed set $C$, and $C,\mathbb{A}^1_2$ must have disjoint projections. If 
$\mathbb{X}^0_2,\mathbb{X}^1_2$ are compact, then $C$ and its projections are compact too. The product of these compact projections is a $(\boratwo\!\times\!\boratwo )_\sigma$ set separating 
$\mathbb{A}^0_2$ from $\mathbb{A}^1_2$, which cannot be. This fact implies that we cannot extend the continuous maps of Theorem 1.5.(e) to $3^\omega$ in general.\bigskip
 
\noindent\bf Notation.\rm\ We now recall some facts about the Gandy-Harringtion topology (see [L2]). Let $Z$ be a recursively presented Polish space. The 
$Gandy$-$Harrington\ topology$ on $Z$ is generated by the $\Ana$ subsets of $Z$. We set  
$\Omega :=\{ z\!\in\! Z\mid\omega_1^z\! =\!\omega_1^{\hbox{\rm CK}}\}$. Then 
$\Omega$ is $\Ana$, dense in $(Z,\mbox{GH})$, and $W\cap\Omega$ is a clopen subset of 
$(\Omega ,{\it\Sigma}_Z )$ for each $W\!\in\!\Ana (Z)$. Moreover, $(\Omega ,\mbox{GH})$ is a $0$-dimensional Polish space. So we fix a complete compatible metric $d_{\mbox{GH}}$ on 
$(\Omega ,\mbox{GH})$.\bigskip

\noindent\bf Proof of Theorem 1.5.\rm\ We already saw that (a)-(d) are equivalent at the end of Section 2. Lemma 4.3 shows that (e) implies (a). So it is enough to show that (d) implies (e). We set 
$N\! :=\! A_0\cap\overline{A_1}^{T_2\times T_2}$, which is not empty. Lemma 2.2 implies that
$$\begin{array}{ll}
(x,y)\!\notin\!\overline{A_1}^{T_2\times T_2}\!\!\!\!
& \Leftrightarrow\exists C,D\!\in\!\Ana\cap\bormone ~~(x,y)\!\in\! C\!\times\! D\!\subseteq\!\neg A_1\cr
& \Leftrightarrow\exists C,D\!\in\!\Borel\cap\boratwo ~~(x,y)\!\in\! C\!\times\! D\!\subseteq\!\neg A_1.
\end{array}$$
This and Proposition 1.4 in [Lo1] show that $N$ is $\Ana$.\bigskip

\noindent $\bullet$ Note that $s$ is not suitable if and only if it is of the form $u\varepsilon v$, where $u$ is suitable, $\varepsilon\!\in\! 2$ and $v\!\in\! 2^{<\omega}$. If $\emptyset\!\not=\! s$ is suitable, then we set $s^-\! :=\! s\vert\mbox{max}\{ l\! <\!\vert s\vert\mid s\vert l\mbox{ is suitable}\}$. We construct\smallskip 

\noindent - a sequence $(x_s)_{s\in 3^{<\omega}}$ of points of $X$,\smallskip

\noindent - a sequence $(y_s)_{s\in 3^{<\omega}}$ of points of $Y$,\smallskip

\noindent - a sequence $(U_s)_{s\in 3^{<\omega}}$ of $\Boraone$ subsets of $X$,\smallskip

\noindent - a sequence $(V_s)_{s\in 3^{<\omega}}$ of $\Boraone$ subsets of $Y$,\smallskip

\noindent - a sequence $(W_s)_{s\in 3^{<\omega}\mbox{ suitable}}$ of $\Ana$ subsets of $X\!\times\!Y$.
\bigskip

 We want these objects to satisfy the following conditions:
$$\begin{array}{ll}
& (1)~(x_s,y_s)\!\in\! U_s\!\times\! V_s\cr
& (2)~(x_s,y_s)\!\in\! W_s\!\subseteq\! N\cap\Omega\mbox{ if }s\mbox{ is suitable}\cr
& (3)~\overline{U_{s\varepsilon}}\!\subseteq\! U_s\mbox{ if }s\mbox{ is suitable or }s\! =\! u0v
\mbox{, and }\overline{U_{u1v2}}\!\subseteq\! U_u\cr
& (4)~\overline{V_{s\varepsilon}}\!\subseteq\! V_s\mbox{ if }s\mbox{ is suitable or }s\! =\! u1v
\mbox{, and }\overline{V_{u0v2}}\!\subseteq\! V_u\cr
& (5)~W_s\!\subseteq\! W_{s^-}\mbox{ if }\emptyset\!\not=\! s\mbox{ is suitable}\cr
& (6)~\mbox{diam}(U_s),\mbox{diam}(V_s)\!\leq\! 2^{-\vert s\vert}\cr
& (7)~\mbox{diam}_{\mbox{GH}}(W_s)\!\leq\! 2^{-\vert s\vert}\mbox{ if }s\mbox{ is suitable}\cr
& (8)~(x_{u0},y_{u1})\!\in\!\big(\overline{\Pi_0[(U_u\!\times\! V_u)\cap W_u]}\!\times\!
\overline{\Pi_1[(U_u\!\times\! V_u)\cap W_u]}\big)\cap A_1\cr
& (9)~(x_{u0v},y_{u1v})\! =\! (x_{u0},y_{u1})
\end{array}$$ 

\vfill\eject

\noindent $\bullet$ Assume that this is done. Let $\alpha\!\in\!\mathbb{X}^0_2$. Then the increasing sequence $(p_k)$ of integers such that $\alpha\vert p_k$ is suitable or of the form $u0v$ is infinite. Condition (3) implies that $(\overline{U_{\alpha\vert p_k}})_{k\in\omega}$ is non-increasing. Moreover, 
$(\overline{U_{\alpha\vert p_k}})_{k\in\omega}$ is a sequence of nonempty closed subsets of $X$ whose diameters tend to $0$, so that we can define 
$\{ f(\alpha )\}\! :=\!\bigcap_{k\in\omega}~\overline{U_{\alpha\vert p_k}}\! =\!
\bigcap_{k\in\omega}~U_{\alpha\vert p_k}$. This defines a continuous map 
$f\! :\!\mathbb{X}^0_2\!\rightarrow\! X$ with 
$f(\alpha )\! =\!\mbox{lim}_{k\rightarrow\infty}~x_{\alpha\vert p_k}$. Similarly, we define 
$g\! :\!\mathbb{X}^1_2\!\rightarrow\! Y$ continuous with 
$g(\beta )\! =\!\mbox{lim}_{k\rightarrow\infty}~y_{\beta\vert q_k}$.\bigskip

 If $\alpha\!\in\!\mathbb{X}^0_2\cap\mathbb{X}^1_2$, then the sequence 
$(k_j)$ of integers such that $\alpha\vert p_{k_j}$ is suitable is infinite. Note that 
$(W_{\alpha\vert p_{k_j}})_{j\in\omega}$ is a non-increasing sequence of nonempty closed subsets of 
$\Omega$ whose GH-diameters tend to $0$, so that we can define $F(\alpha )$ by 
$\{ F(\alpha )\}\! :=\!\bigcap_{j\in\omega}~W_{\alpha\vert p_{k_j}}\!\subseteq\! N\!\subseteq\! A_0$. As 
$F(\alpha )$ is the limit (in $(X\!\times\! Y,\mbox{GH})$, and thus in $X\!\times\! Y$) of 
$(x_{\alpha\vert p_{k_j}},y_{\alpha\vert p_{k_j}})_{j\in\omega}$, we get 
$F(\alpha )\! =\!\big( f(\alpha ),g(\alpha )\big)$. Thus 
$\mathbb{A}^0_2\!\subseteq\! (f\!\times\!g)^{-1}(A_0)$.\bigskip

 Note that $x_{s\varepsilon}\! =\! x_{s\varepsilon (\alpha\vert 1)}\! =\! ...\! =\! 
x_{s\varepsilon (\alpha\vert (q+1))}$ for each 
$(s,\varepsilon ,\alpha )\!\in\! 3^{<\omega}\!\times\! 2\!\times\! 2^\omega$. This implies that  
$f(s0\alpha )\! =\!\mbox{lim}_{q\rightarrow\infty}~x_{s0 (\alpha\vert q)}\! =\! x_{s0}$. Similarly, 
$g(s1\beta )\! =\! y_{s1}$ and $\big( f(s0\alpha ),g(s1\beta )\big)\! =\! (x_{s0},y_{s1})\!\in\! A_1$. 
Thus $\mathbb{A}^1_2\!\subseteq\! (f\!\times\!g)^{-1}(A_1)$.\bigskip

\noindent $\bullet$ Let us prove that the construction is possible. As $N$ is not empty, we can choose 
$(x_\emptyset ,y_\emptyset )\!\in\! N\cap\Omega$, a $\Ana$ subset $W_\emptyset$ of $X\!\times\! Y$ with 
$(x_\emptyset ,y_\emptyset )\!\in\! W_\emptyset\!\subseteq\! N\cap\Omega$ of GH-diameter at most $1$, and a $\Boraone$ neighborhood $U_\emptyset$ (resp., $V_\emptyset$) of $x_\emptyset$ (resp., $y_\emptyset$) of diameter at most $1$. Assume that $(x_s)_{s\in 3^{\leq l}}$, 
$(y_s)_{s\in 3^{\leq l}}$, $(U_s)_{s\in 3^{\leq l}}$, $(V_s)_{s\in 3^{\leq l}}$ and $(W_s)_{s\in 3^{\leq l}}$ satisfying (1)-(9) have been constructed, which is the case for $l\! =\! 0$.\bigskip

 Note that $(x_u,y_u)\!\in\! (U_u\!\times\! V_u)\cap W_u\!\subseteq\! 
\overline{A_1}^{T_2\times T_2}$ since $u$ is suitable. We choose $U,V\!\in\!\Boraone$ with 
$(x_u,y_u)\!\in\! U\!\times\! V\!\subseteq\!\overline{U}\!\times\!\overline{V}\!\subseteq\! U_u\!\times\! V_u$. As $\Pi_\varepsilon [(U\!\times\! V)\cap W_u]$ is $\Ana$, 
$\overline{\Pi_\varepsilon [(U\!\times\! V)\cap W_u]}$ is $\Ana\cap\bormone$. In particular, $\overline{\Pi_\varepsilon [(U\!\times\! V)\cap W_u]}$ is $T_2$-open. This shows the existence of $$(x_{u0},y_{u1})\!\in\!\big(\overline{\Pi_0[(U\!\times\! V)\cap W_u]}\!\times\!
\overline{\Pi_1[(U\!\times\! V)\cap W_u]}\big)\cap A_1.$$ 
Note that $(x_{u0},y_{u1})\!\in\!\overline{U}\!\times\!\overline{V}\!\subseteq\! U_u\!\times\! V_u$. We set 
$x_{u1}\! :=\! x_u$, $y_{u0}\! :=\! y_u$. We defined $x_s,y_s$ when $s\!\in\! 3^{l+1}$ is not suitable but $s\vert l$ is suitable.\bigskip

 Assume now that $s$ is suitable, but not $s\vert l$. This gives $(u,\varepsilon ,v)$ such that 
$s\! =\! u\varepsilon v2$. Assume first that $\varepsilon\! =\! 0$. Note that 
$x_{u0v}\! =\! x_{u0}\!\in\! U_{u0v}\cap\overline{\Pi_0[(U_u\!\times\! V_u)\cap W_u]}$. This gives 
$x_s\!\in\! U_{u0v}\cap\Pi_0[(U_u\!\times\! V_u)\cap W_u]$, and also $y_s$ with 
$(x_s,y_s)\!\in\!\big( (U_u\cap U_{u0v})\!\times\! V_u\big)\cap W_u\! =\! 
(U_{u0v}\!\times\! V_u)\cap W_u$. If $\varepsilon\! =\! 1$, then similarly we get 
$(x_s,y_s)\!\in\! (U_u\!\times\! V_{u1v})\cap W_u$.\bigskip

 If $s$ and $s\vert l$ are both suitable, or both non suitable, then we set 
$(x_s,y_s)\! :=\! (x_{s\vert l},y_{s\vert l})$. So we defined $x_s,y_s$ in any case. Note that Conditions 
(8) and (9) are fullfilled, and that $(x_s,y_s)\!\in\! W_{s^-}$ if $s$ is suitable. Moreover, 
$x_s\!\in\! U_{s\vert l}$ if $s\vert l$ is suitable or $s\vert l\! =\! u0v$, and $x_s\!\in\! U_u$ if 
$s\! =\! u1v2$, and similarly in $Y$. We choose $\Boraone$ sets $U_s,V_s$ of diameter at most 
$2^{-l-1}$ with
$$(x_s,y_s)\!\in\! U_s\!\times\! V_s\!\subseteq\!\overline{U_s}\!\times\!\overline{V_s}\!\subseteq\!
\left\{\!\!\!\!\!\!
\begin{array}{ll}
& U_{s\vert l}\!\times\! V_{s\vert l}\mbox{ if }s\mbox{ is not suitable or }s\vert l\mbox{ is suitable,}\cr
& U_{s\vert l}\!\times\! V_u\mbox{ if }s\! =\! u0v2\mbox{,}\cr
& U_u\!\times\! V_{s\vert l}\mbox{ if }s\! =\! u1v2.
\end{array}
\right.$$
It remains to choose, when $s$ is suitable, $W_s\!\in\!\Ana (X\!\times\! Y)$ of GH-diameter at most 
$2^{-l-1}$ with $(x_s,y_s)\!\in\! W_s\!\subseteq\! W_{s^-}$.\hfill{$\square$}

\vfill\eject
 
\section{$\!\!\!\!\!\!$ References}

\noindent [K]\ \ A. S. Kechris,~\it Classical Descriptive Set Theory,~\rm 
Springer-Verlag, 1995

\noindent [K-S-T]\ \ A. S. Kechris, S. Solecki and S. Todor\v cevi\'c, Borel chromatic numbers,\ \it 
Adv. Math.\rm\ 141 (1999), 1-44

\noindent [L1]\ \ D. Lecomte, On minimal non potentially closed subsets of the plane,
\ \it Topology Appl.~\rm 154, 1 (2007) 241-262

\noindent [L2]\ \ D. Lecomte, A dichotomy characterizing analytic graphs of uncountable Borel chromatic number in any dimension,~\it Trans. Amer. Math. Soc.\rm~361 (2009), 4181-4193

\noindent [L3]\ \ D. Lecomte, How can we recognize potentially $\bormxi$ subsets of the plane?,
~\it to appear in J. Math. Log. (see arXiv)\rm

\noindent [Lo1]\ \ A. Louveau, A separation theorem for $\Ana$ sets,\ \it Trans. 
Amer. Math. Soc.\ \rm 260 (1980), 363-378

\noindent [Lo2]\ \ A. Louveau, Ensembles analytiques et bor\'eliens dans les 
espaces produit,~\it Ast\'erisque (S. M. F.)\ \rm 78 (1980)

\noindent [M]\ \ Y. N. Moschovakis,~\it Descriptive set theory,~\rm North-Holland, 1980\bigskip\bigskip

\noindent\bf Acknowledgements.\rm\ This work started last summer, when I was invited at the University 
of Prague by Miroslav Zelen\'y. I am very grateful to him for that, and pleased that we could prove Theorem 4.2 together.
 
\end{document}